\documentclass[12pt,twoside]{amsart}
\usepackage[latin1]{inputenc}
\usepackage{amsmath, amsthm, amscd, amsfonts, amssymb, graphicx}
\usepackage[bookmarksnumbered, plainpages]{hyperref}

\newtheorem{thm}{Theorem}[section]

\newtheorem{lem}[thm]{Lemma}

\newtheorem{defn}[thm]{Definition}

\numberwithin{equation}{section}

\begin{document}

\title{\bf $\Gamma^{0}(2)$ modular forms and anomaly cancellation formulas}
\author{Siyao Liu \hskip 0.4 true cm  Yong Wang$^{*}$}

\thanks{{\scriptsize
\hskip -0.4 true cm \textit{2010 Mathematics Subject Classification:}
58C20; 57R20; 53C80.
\newline \textit{$^{*}$Corresponding author}}}

\maketitle

\begin{abstract}
 \indent In \cite{HH}, \cite{W1} and \cite{LW2}, the authors gave some modular forms over $\Gamma^0(2).$ In this note, we proceed with the study of cancellation formulas relating to the modular forms.

\end{abstract}

\vskip 0.2 true cm


\pagestyle{myheadings}
\markboth{\rightline {\scriptsize Liu}}
         {\leftline{\scriptsize $\Gamma^{0}(2)$ modular forms and anomaly cancellation formulas}}

\bigskip
\bigskip


\section{ Introduction }

In \cite{AW}, Alvarez-Gaum\'{e} and Witten discovered the ``miraculous cancellation" formula for gravitational anomaly  as follow:
\begin{align}
\Big\{\widehat{L}(TM, \nabla^{TM})\Big\}^{(12)}=\Big\{\widehat{A}(TM, \nabla^{TM})ch(T_{\mathbf{C}}M, \nabla^{T_{\mathbf{C}}M})-32\widehat{A}(TM, \nabla^{TM})\Big\}^{(12)},\nonumber
\end{align}
where $T_{\mathbf{C}}M$ denotes the complexification of $TM$ and $\nabla^{T_{\mathbf{C}}M}$ is canonically induced from $\nabla^{TM},$ the Levi-Civita connection associated to the Riemannian structure of $M.$
This formula reveals a beautiful relation between the top components of the Hirzebruch $\widehat{L}$-form and $\widehat{A}$-form of a $12$-dimensional smooth Riemannian manifold $M.$
Liu established higher-dimensional ``miraculous cancellation" formulas for $(8k+4)$-dimensional Riemannian manifolds by developing modular invariance properties of characteristic forms \cite{L}.
Han and Zhang established a general cancellation formula that involves a complex line bundle for $(8k+4)$-dimensional smooth Riemannian manifold in \cite{HZ1,HZ2}.
For higher-dimensional smooth Riemannian manifolds the authors obtained some cancellation formulas in \cite{HH}.
In \cite{W1}, Wang proved more general cancellation formulas for $(8k+2)$ and $(8k+6)$-dimensional smooth Riemannian manifolds.
And in \cite{LW1}, the authors generalized the Han-Liu-Zhang cancellation formulas to the $(a, b)$ type cancellation formulas.
To enrich the results, Liu and Wang proved some other $(a, b)$ type cancellation formulas for even-dimensional Riemannian manifolds \cite{LW2}.

In this paper, we are interested in finding some new cancellation formulas.
According to \cite{HH,W1} and \cite{LW2}, we found some $\Gamma^0(2)$ modular forms.
Different from the methods used in the articles, we conclude some new cancellation formulas.

A brief description of the organization of this paper is as follows.
In Section 2, we give some definitions and basic notions that we will use in this paper.
In the next section, based on the modular forms over $\Gamma^0(2)$ in \cite{HH,W1}, we compute some cancellation formulas.
Finally, in Section 4, we prove some cancellation formulas and cancellation formulas involving a complex line bundle for $4d$-dimensional Riemannian manifolds.


\vskip 1 true cm

\section{ Characteristic forms and modular forms }

Firstly, we give some definitions and basic notions on characteristic forms and modular forms that will be used throughout the paper.
For the details, see \cite{A,H,Z}.

\vskip 0.3 true cm
2.1. \noindent{ Characteristic forms }
\vskip 0.3 true cm

Let $M$ be a Riemannian manifold, $\nabla^{TM}$ be the associated Levi-Civita connection on $TM$ and $R^{TM}=(\nabla^{TM})^{2}$ be the curvature of $\nabla^{TM}.$
According to the detailed descriptions in \cite{Z}, let $\widehat{A}(TM, \nabla^{TM})$ and $\widehat{L}(TM, \nabla^{TM})$ be the Hirzebruch characteristic forms defined respectively by
\begin{align}
\widehat{A}(TM, \nabla^{TM})&=\det\nolimits^{\frac{1}{2}}\left(\frac{\frac{\sqrt{-1}}{4\pi}R^{TM}}{\sinh(\frac{\sqrt{-1}}{4\pi}R^{TM})}\right),\\
\widehat{L}(TM, \nabla^{TM})&=\det\nolimits^{\frac{1}{2}}\left(\frac{\frac{\sqrt{-1}}{2\pi}R^{TM}}{\tanh(\frac{\sqrt{-1}}{4\pi}R^{TM})}\right).
\end{align}

Let $E, F$ be two Hermitian vector bundles over $M$ carrying Hermitian connection $\nabla^{E}, \nabla^{F}$ respectively.
Let $R^{E}=(\nabla^{E})^{2}$ (resp. $R^{F}=(\nabla^{F})^{2}$) be the curvature of $\nabla^{E}$ (resp. $\nabla^{F}$).
If we set the formal difference $G=E-F,$ then $G$ carries an induced Hermitian connection $\nabla^{G}$ in an obvious sense.
We define the associated Chern character form as
\begin{align}
\text{ch}(G, \nabla^{G})=\text{tr}\left[\exp\left(\frac{\sqrt{-1}}{2\pi}R^{E}\right)\right]-\text{tr}\left[\exp\left(\frac{\sqrt{-1}}{2\pi}R^{F}\right)\right].
\end{align}

For any complex number $t,$ let
\begin{align}
\wedge_{t}(E)&=\mathbf{C}|_{M}+tE+t^{2}\wedge^{2}(E)+\cdot\cdot\cdot,\\
S_{t}(E)&=\mathbf{C}|_{M}+tE+t^{2}S^{2}(E)+\cdot\cdot\cdot,
\end{align}
denote respectively the total exterior and symmetric powers of $E,$ which live in $K(M)[[t]].$
The following relations between these operations hold,
\begin{align}
S_{t}(E)=\frac{1}{\wedge_{-t}(E)},~~~~ \wedge_{t}(E-F)=\frac{\wedge_{t}(E)}{\wedge_{t}(F)}.
\end{align}
Moreover, if $\{ \omega_{i} \}, \{ \omega'_{j} \}$ are formal Chern roots for Hermitian vector bundles $E, F$ respectively, then
\begin{align}
\text{ch}(\wedge_{t}(E))=\prod_{i}(1+e^{\omega_{i}}t).
\end{align}
Then we have the following formulas for Chern character forms,
\begin{align}
\text{ch}(S_{t}(E))=\frac{1}{\prod\limits_{i}(1-e^{\omega_{i}}t)},~~~ \text{ch}(\wedge_{t}(E-F))=\frac{\prod\limits_{i}(1+e^{\omega_{i}}t)}{\prod\limits_{j}(1+e^{\omega'_{j}}t)}.
\end{align}

If $W$ is a real Euclidean vector bundle over $M$ carrying a Euclidean connection $\nabla^{W},$ then its complexification $W_{\mathbf{C}}=W\otimes \mathbf{C}$ is a complex vector bundle over $M$ carrying a canonical induced Hermitian metric from that of $W,$ as well as a Hermitian connection $\nabla^{W_{\mathbf{C}}}$ induced from $\nabla^{W}.$
If $E$ is a vector bundle (complex or real) over $M,$ set $\widetilde{E}=E-\dim E$ in $K(M)$ or $KO(M).$

\vskip 0.3 true cm
2.2. \noindent{ Some properties about the Jacobi theta functions and modular forms }
\vskip 0.3 true cm

Refer to \cite{C}, we recall the four Jacobi theta functions are defined as follows:
\begin{align}
&\theta(v, \tau)=2q^{\frac{1}{8}}\sin(\pi v)\prod^{\infty}_{j=1}[(1-q^{j})(1-e^{2\pi\sqrt{-1}v}q^{j})(1-e^{-2\pi\sqrt{-1}v}q^{j})];\\
&\theta_{1}(v, \tau)=2q^{\frac{1}{8}}\cos(\pi v)\prod^{\infty}_{j=1}[(1-q^{j})(1+e^{2\pi\sqrt{-1}v}q^{j})(1+e^{-2\pi\sqrt{-1}v}q^{j})];\\
&\theta_{2}(v, \tau)=\prod^{\infty}_{j=1}[(1-q^{j})(1-e^{2\pi\sqrt{-1}v}q^{j-\frac{1}{2}})(1-e^{-2\pi\sqrt{-1}v}q^{j-\frac{1}{2}})];\\
&\theta_{3}(v, \tau)=\prod^{\infty}_{j=1}[(1-q^{j})(1+e^{2\pi\sqrt{-1}v}q^{j-\frac{1}{2}})(1+e^{-2\pi\sqrt{-1}v}q^{j-\frac{1}{2}})],
\end{align}
where $q=e^{2\pi\sqrt{-1}\tau}$ with $\tau\in\mathbf{H},$ the upper half complex plane.

Let
\begin{align}
\theta'(0, \tau)=\frac{\partial \theta(v, \tau)}{\partial v}\Big|_{v=0}.
\end{align}
Then the following Jacobi identity holds,
\begin{align}
\theta'(0, \tau)=\pi\theta_{1}(0, \tau)\theta_{2}(0, \tau)\theta_{3}(0, \tau).
\end{align}

In what follows,
\begin{align}
SL_{2}(\mathbf{Z})=\Big\{ \Big( \begin{array}{cc}
a & b\\
c & d
\end{array}
\Big)\Big|a, b, c, d\in\mathbf{Z}, ad-bc=1
\Big \}
\end{align}
stands for the modular group.
Write $S=\Big(\begin{array}{cc}
0 & -1\\
1 & 0
\end{array}\Big),
T=\Big(\begin{array}{cc}
1 & 1\\
0 & 1
\end{array}\Big)$
be the two generators of $SL_{2}(\mathbf{Z}).$
They act on $\mathbf{H}$ by $S\tau=-\frac{1}{\tau}, T\tau=\tau+1.$
One has the following transformation laws of theta functions under the actions of $S$ and $T:$
\begin{align}
&\theta(v, \tau+1)=e^{\frac{\pi\sqrt{-1}}{4}}\theta(v, \tau),~~~~ \theta(v, -\frac{1}{\tau})=\frac{1}{\sqrt{-1}}\Big(\frac{\tau}{\sqrt{-1}}\Big)^{\frac{1}{2}}e^{\pi\sqrt{-1}\tau v^{2}}\theta(\tau v, \tau);\\
&\theta_{1}(v, \tau+1)=e^{\frac{\pi\sqrt{-1}}{4}}\theta_{1}(v, \tau),~~~~ \theta_{1}(v, -\frac{1}{\tau})=\Big(\frac{\tau}{\sqrt{-1}}\Big)^{\frac{1}{2}}e^{\pi\sqrt{-1}\tau v^{2}}\theta_{2}(\tau v, \tau);\\
&\theta_{2}(v, \tau+1)=\theta_{3}(v, \tau),~~~~ \theta_{2}(v, -\frac{1}{\tau})=\Big(\frac{\tau}{\sqrt{-1}}\Big)^{\frac{1}{2}}e^{\pi\sqrt{-1}\tau v^{2}}\theta_{1}(\tau v, \tau);\\
&\theta_{3}(v, \tau+1)=\theta_{2}(v, \tau),~~~~ \theta_{3}(v, -\frac{1}{\tau})=\Big(\frac{\tau}{\sqrt{-1}}\Big)^{\frac{1}{2}}e^{\pi\sqrt{-1}\tau v^{2}}\theta_{3}(\tau v, \tau).
\end{align}
Differentiating the above transformation formulas, we get that
\begin{align}
&\theta'(v, \tau+1)=e^{\frac{\pi\sqrt{-1}}{4}}\theta'(v, \tau),\\
&\theta'(v, -\frac{1}{\tau})=\frac{1}{\sqrt{-1}}\Big(\frac{\tau}{\sqrt{-1}}\Big)^{\frac{1}{2}}e^{\pi\sqrt{-1}\tau v^{2}}(2\pi\sqrt{-1}\tau v\theta(\tau v, \tau)+\tau\theta'(\tau v, \tau));\nonumber\\
&\theta_{1}'(v, \tau+1)=e^{\frac{\pi\sqrt{-1}}{4}}\theta_{1}'(v, \tau),\\
&\theta_{1}'(v, -\frac{1}{\tau})=\Big(\frac{\tau}{\sqrt{-1}}\Big)^{\frac{1}{2}}e^{\pi\sqrt{-1}\tau v^{2}}(2\pi\sqrt{-1}\tau v\theta_{2}(\tau v, \tau)+\tau\theta'_{2}(\tau v, \tau));\nonumber\\
&\theta'_{2}(v, \tau+1)=\theta'_{3}(v, \tau),\\
&\theta'_{2}(v, -\frac{1}{\tau})=\Big(\frac{\tau}{\sqrt{-1}}\Big)^{\frac{1}{2}}e^{\pi\sqrt{-1}\tau v^{2}}(2\pi\sqrt{-1}\tau v\theta_{1}(\tau v, \tau)+\tau\theta'_{1}(\tau v, \tau));\nonumber\\
&\theta'_{3}(v, \tau+1)=\theta'_{2}(v, \tau),\\
&\theta'_{3}(v, -\frac{1}{\tau})=\Big(\frac{\tau}{\sqrt{-1}}\Big)^{\frac{1}{2}}e^{\pi\sqrt{-1}\tau v^{2}}(2\pi\sqrt{-1}\tau v\theta_{3}(\tau v, \tau)+\tau\theta'_{3}(\tau v, \tau)).\nonumber
\end{align}
Therefore
\begin{align}
\theta'(0, -\frac{1}{\tau})=\frac{1}{\sqrt{-1}}\Big(\frac{\tau}{\sqrt{-1}}\Big)^{\frac{1}{2}}\tau\theta'(0, \tau).
\end{align}

\begin{defn}
A modular form over $\Gamma,$ a subgroup of $SL_{2}(\mathbf{Z}),$ is a holomorphic function $f(\tau)$ on $\mathbf{H}$ such that
\begin{align}
f(g\tau):=f\Big(\frac{a\tau+b}{c\tau+d}\Big)=\chi(g)(c\tau+d)^{k}f(\tau),~\forall g=\Big(\begin{array}{cc}
a & b\\
c & d
\end{array}\Big)\in\Gamma,
\end{align}
where $\chi:\Gamma\rightarrow\mathbf{C}^{*}$ is a character of $\Gamma,$ $k$ is called the weight of $f.$
\end{defn}

Let
\begin{align}
&\Gamma_{0}(2)=\Big\{ \Big( \begin{array}{cc}
a & b\\
c & d
\end{array}
\Big)\in SL_{2}(\mathbf{Z})\Big|c\equiv0(\text{mod}~2)
\Big \},\\
&\Gamma^{0}(2)=\Big\{ \Big( \begin{array}{cc}
a & b\\
c & d
\end{array}
\Big)\in SL_{2}(\mathbf{Z})\Big|b\equiv0(\text{mod}~2)
\Big \},\\
&\Gamma_{\theta}=\Big\{ \Big( \begin{array}{cc}
a & b\\
c & d
\end{array}
\Big)\in SL_{2}(\mathbf{Z})\Big|\Big( \begin{array}{cc}
a & b\\
c & d
\end{array}
\Big)\equiv\Big( \begin{array}{cc}
1 & 0\\
0 & 1
\end{array}
\Big) \text{or} \Big( \begin{array}{cc}
0 & 1\\
1 & 0
\end{array}
\Big)(\text{mod}~2)
\Big \}
\end{align}
be the three modular subgroups of $SL_{2}(\mathbf{Z}).$
It is known that the generators of $\Gamma_{0}(2)$ are $T, ST^{2}ST,$ the generators of $\Gamma^{0}(2)$ are $STS, T^{2}STS$ and the generators of $\Gamma_{\theta}$ are $ S, T^{2}.$

Let $E_{2}(\tau)$ be Eisenstein series which is a quasimodular form over $SL_{2}(\mathbf{Z}),$ satisfying
\begin{align}
E_{2}\Big(\frac{a\tau+b}{c\tau+d}\Big)=(c\tau+d)^{2}E_{2}(\tau)-\frac{6\sqrt{-1}c(c\tau+d)}{\pi}.
\end{align}
In particular, we have
\begin{align}
&E_{2}(\tau+1)=E_{2}(\tau),\\
&E_{2}\Big(-\frac{1}{\tau}\Big)=\tau^{2}E_{2}(\tau)-\frac{6\sqrt{-1}\tau}{\pi},
\end{align}
and
\begin{align}
E_{2}(\tau)=1-24q-72q^{2}+\cdot\cdot\cdot,
\end{align}
where the $``\cdot\cdot\cdot"$ terms are the higher degree terms, all of which have integral coefficients.

If $\Gamma$ is a modular subgroup, let $\mathcal{M}_{\mathbf{R}}(\Gamma)$ denote the ring of modular forms over $\Gamma$ with real Fourier coefficients.
Writing $\theta_{j}=\theta_{j}(0,\tau), 1\leq j\leq3,$ we introduce six explicit modular forms,
\begin{align}
&\delta_{1}(\tau)=\frac{1}{8}(\theta^{4}_{2}+\theta^{4}_{3}),~~~~\varepsilon_{1}(\tau)=\frac{1}{16}\theta^{4}_{2}\theta^{4}_{3};\\
&\delta_{2}(\tau)=-\frac{1}{8}(\theta^{4}_{1}+\theta^{4}_{3}),~~~~\varepsilon_{2}(\tau)=\frac{1}{16}\theta^{4}_{1}\theta^{4}_{3};\\
&\delta_{3}(\tau)=\frac{1}{8}(\theta^{4}_{1}-\theta^{4}_{2}),~~~~\varepsilon_{3}(\tau)=-\frac{1}{16}\theta^{4}_{1}\theta^{4}_{2}.
\end{align}
They have the following Fourier expansions in $q^{\frac{1}{2}}$:
\begin{align}
&\delta_{1}(\tau)=\frac{1}{4}+6q+\cdot\cdot\cdot,~~~~\varepsilon_{1}(\tau)=\frac{1}{16}-q+\cdot\cdot\cdot;\\
&\delta_{2}(\tau)=-\frac{1}{8}-3q^{\frac{1}{2}}+\cdot\cdot\cdot,~~~~\varepsilon_{2}(\tau)=q^{\frac{1}{2}}+\cdot\cdot\cdot;\\
&\delta_{3}(\tau)=-\frac{1}{8}+3q^{\frac{1}{2}}+\cdot\cdot\cdot,~~~~\varepsilon_{3}(\tau)=-q^{\frac{1}{2}}+\cdot\cdot\cdot.
\end{align}
They also satisfy the transformation laws,
\begin{align}
&\delta_{2}\Big(-\frac{1}{\tau}\Big)=\tau^{2}\delta_{1}(\tau),~~~~\varepsilon_{2}\Big(-\frac{1}{\tau}\Big)=\tau^{4}\varepsilon_{1}(\tau);\\
&\delta_{2}(\tau+1)=\delta_{3}(\tau),~~~~\varepsilon_{2}(\tau+1)=\varepsilon_{3}(\tau).
\end{align}

\begin{lem}(\cite{L})
$\delta_{1}(\tau)$ (resp. $\varepsilon_{1}(\tau)$) is a modular form of weight 2 (resp. 4) over $\Gamma_{0}(2),$ $\delta_{2}(\tau)$ (resp. $\varepsilon_{2}(\tau)$) is a modular form of weight 2 (resp. 4) over $\Gamma^{0}(2),$ while $\delta_{3}(\tau)$ (resp. $\varepsilon_{3}(\tau)$) is a modular form of weight 2 (resp. 4) over $\Gamma_{\theta},$ and moreover $\mathcal{M}_{\mathbf{R}}(\Gamma^{0}(2))=\mathbf{R}[\delta_{2}(\tau), \varepsilon_{2}(\tau)].$
\end{lem}

\section{ Cancellation formulas }

\vskip 0.3 true cm
3.1. \noindent{ Cancellation formulas for $8k+2$-dimensional Riemannian manifolds  }
\vskip 0.3 true cm

Let $M$ be a $8k+2$-dimensional Riemannian manifold and $\xi$ be a rank two real oriented Euclidean vector bundle over $M$ carrying with a Euclidean connection $\nabla^\xi.$

Define
\begin{align}
\Theta_{2}(T_{\mathbf{C}}M+\xi_{\mathbf{C}}, \xi_{\mathbf{C}})=&\bigotimes^{\infty}_{n=1}S_{q^{n}}(\widetilde{T_{\mathbf{C}}M}+\widetilde{\xi_{\mathbf{C}}})\otimes
\bigotimes_{m=1}^{\infty}\wedge_{-q^{m-\frac{1}{2}}}(\widetilde{T_{\mathbf{C}}M}+\widetilde{\xi_{\mathbf{C}}}-2\widetilde{\xi_{\mathbf{C}}})\\
&\otimes \bigotimes _{r=1}^{\infty}\wedge_{q^{r-\frac{1}{2}}}(\widetilde{\xi_{\mathbf{C}}})\otimes\bigotimes_{s=1}^{\infty}\wedge _{q^{s}}(\widetilde{\xi_{\mathbf{C}}}),\nonumber\\
\Theta_{2}(T_{\mathbf{C}}M+\xi_{\mathbf{C}}, \mathbf{C}^{2})=&\bigotimes^{\infty}_{n=1}S_{q^{n}}(\widetilde{T_{\mathbf{C}}M}+\widetilde{\xi_{\mathbf{C}}})\otimes
\bigotimes_{m=1}^{\infty}\wedge_{-q^{m-\frac{1}{2}}}(\widetilde{T_{\mathbf{C}}M}+\widetilde{\xi_{\mathbf{C}}}).
\end{align}
$\Theta_{2}(T_{\mathbf{C}}M+\xi_{\mathbf{C}}, \xi_{\mathbf{C}})$ and $\Theta_{2}(T_{\mathbf{C}}M+\xi_{\mathbf{C}}, \mathbf{C}^{2})$ admit formal Fourier expansion in $q^{\frac{1}{2}}$ as
\begin{align}
\Theta_{2}(T_{\mathbf{C}}M+\xi_{\mathbf{C}}, \xi_{\mathbf{C}})=&A_{0}(T_{\mathbf{C}}M+\xi_{\mathbf{C}}, \xi_{\mathbf{C}})+A_{1}(T_{\mathbf{C}}M+\xi_{\mathbf{C}}, \xi_{\mathbf{C}})q^{\frac{1}{2}}+\cdot\cdot\cdot,\\
\Theta_{2}(T_{\mathbf{C}}M+\xi_{\mathbf{C}}, \mathbf{C}^{2})=&A_{0}(T_{\mathbf{C}}M+\xi_{\mathbf{C}}, \mathbf{C}^{2})+A_{1}(T_{\mathbf{C}}M+\xi_{\mathbf{C}}, \mathbf{C}^{2})q^{\frac{1}{2}}+\cdot\cdot\cdot,
\end{align}
where $A_{j}$ are elements in the semi-group formally generated by Hermitian vector bundles over $M.$

Set
\begin{align}
P_{2}(\tau)=&\bigg\{\widehat{A}(TM,\nabla^{TM})\frac{1}{2\sinh\Big(\frac{1}{2}c\Big)}ch(\Theta_{2}(T_{\mathbf{C}}M+\xi_{\mathbf{C}},\mathbf{C}^{2}))\\
&-\widehat{A}(TM,\nabla^{TM})\frac{\cosh\Big(\frac{1}{2}c\Big)}{2\sinh\Big(\frac{1}{2}c\Big)}ch(\Theta_2(T_{\mathbf{C}}M+\xi_{\mathbf{C}},\xi_{\mathbf{C}}))\bigg\}^{(8k+2)},\nonumber
\end{align}
where $c=2\pi\sqrt{-1}u.$
Refer to \cite{HH,W1}, we have $P_{2}(\tau)$ is a modular form of weight $4k+2$ over $\Gamma^{0}(2).$

When $k=1,$ by Lemma 2.2, there exist $h_{0}, h_{1}$ such that
\begin{align}
P_{2}(\tau)=&h_{0}(8\delta_{2}(\tau))^{3}+h_{1}(8\delta_{2}(\tau))\varepsilon_{2}(\tau)\\
=&h_{0}(-1-24q^{\frac{1}{2}}-24q+\cdot\cdot\cdot)^{3}+h_{1}(-1-24q^{\frac{1}{2}}-24q+\cdot\cdot\cdot)(q^{\frac{1}{2}}+8q+\cdot\cdot\cdot).\nonumber
\end{align}
Comparing the coefficients of $1, q^{\frac{1}{2}}$ and $q$ in both sides of (3.6), we can get the following formula.
\begin{thm}
For $10$-dimensional Riemannian manifold, we have
\begin{align}
&-504\bigg\{\widehat{A}(TM,\nabla^{TM})\frac{1-\cosh\Big(\frac{1}{2}c\Big)}{2\sinh\Big(\frac{1}{2}c\Big)}\bigg\}^{(10)}+32\bigg\{\widehat{A}(TM,\nabla^{TM})\frac{1}{2\sinh\Big(\frac{1}{2}c\Big)}(ch(-T_{\mathbf{C}}M-\xi_{\mathbf{C}}\\
&+12)-\cosh\Big(\frac{1}{2}c\Big)ch(-T_{\mathbf{C}}M+2\xi_{\mathbf{C}}+6))\bigg\}^{(10)}=\bigg\{\widehat{A}(TM,\nabla^{TM})\frac{1}{2\sinh\Big(\frac{1}{2}c\Big)}(ch(\wedge^{2}(T_{\mathbf{C}}M)\nonumber\\
&+\wedge^{2}(\xi_{\mathbf{C}})+T_{\mathbf{C}}M\otimes\xi_{\mathbf{C}}-11T_{\mathbf{C}}M-11\xi_{\mathbf{C}}+66)-\cosh\Big(\frac{1}{2}c\Big)ch(\wedge^{2}(T_{\mathbf{C}}M)-2T_{\mathbf{C}}M\otimes\xi_{\mathbf{C}}\nonumber\\
&+2\xi_{\mathbf{C}}\otimes\xi_{\mathbf{C}}-5T_{\mathbf{C}}M+14\xi_{\mathbf{C}}+9))\bigg\}^{(10)}.\nonumber
\end{align}
\end{thm}

\vskip 0.3 true cm
3.2. \noindent{ Cancellation formulas for $8k+6$-dimensional Riemannian manifolds  }
\vskip 0.3 true cm

Let $M$ be a $8k+6$-dimensional Riemannian manifold and $\xi$ be a rank two real oriented Euclidean vector bundle over $M$ carrying with a Euclidean connection $\nabla^\xi.$

Now
\begin{align}
\Theta_{2}(T_{\mathbf{C}}M+\xi_{\mathbf{C}}, \xi_{\mathbf{C}})=&A'_{0}(T_{\mathbf{C}}M+\xi_{\mathbf{C}}, \xi_{\mathbf{C}})+A'_{1}(T_{\mathbf{C}}M+\xi_{\mathbf{C}}, \xi_{\mathbf{C}})q^{\frac{1}{2}}+\cdot\cdot\cdot,\\
\Theta_{2}(T_{\mathbf{C}}M+\xi_{\mathbf{C}}, \mathbf{C}^{2})=&A'_{0}(T_{\mathbf{C}}M+\xi_{\mathbf{C}}, \mathbf{C}^{2})+A'_{1}(T_{\mathbf{C}}M+\xi_{\mathbf{C}}, \mathbf{C}^{2})q^{\frac{1}{2}}+\cdot\cdot\cdot,
\end{align}
where $A'_{j}$ are elements in the semi-group formally generated by Hermitian vector bundles over $M.$

Fix
\begin{align}
P'_{2}(\tau)=&\bigg\{\widehat{A}(TM,\nabla^{TM})\frac{1}{2\sinh\Big(\frac{1}{2}c\Big)}ch(\Theta_{2}(T_{\mathbf{C}}M+\xi_{\mathbf{C}},\mathbf{C}^{2}))\\
&-\widehat{A}(TM,\nabla^{TM})\frac{\cosh\Big(\frac{1}{2}c\Big)}{2\sinh\Big(\frac{1}{2}c\Big)}ch(\Theta_2(T_{\mathbf{C}}M+\xi_{\mathbf{C}},\xi_{\mathbf{C}}))\bigg\}^{(8k+6)}.\nonumber
\end{align}
From \cite{HH,W1}, we see that $P'_{2}(\tau)$ is a modular form of weight $4k+4$ over $\Gamma^{0}(2).$

Choose $k=1,$ we obtain
\begin{align}
P'_{2}(\tau)=&h_{0}(8\delta_{2}(\tau))^{4}+h_{1}(8\delta_{2}(\tau))^{2}\varepsilon_{2}(\tau)+h_{2}(\varepsilon_{2}(\tau))^{2}.
\end{align}
Then similar to Theorem 3.1, the following equation holds.
\begin{thm}
For $14$-dimensional Riemannian manifold, we conclude that
\begin{align}
&-7680\bigg\{\widehat{A}(TM,\nabla^{TM})\frac{1-\cosh\Big(\frac{1}{2}c\Big)}{2\sinh\Big(\frac{1}{2}c\Big)}\bigg\}^{(14)}+140\bigg\{\widehat{A}(TM,\nabla^{TM})\frac{1}{2\sinh\Big(\frac{1}{2}c\Big)}(ch(-T_{\mathbf{C}}M\\
&-\xi_{\mathbf{C}}+16)-\cosh\Big(\frac{1}{2}c\Big)ch(-T_{\mathbf{C}}M+2\xi_{\mathbf{C}}+10))\bigg\}^{(14)}+16\bigg\{\widehat{A}(TM,\nabla^{TM})\frac{1}{2\sinh\Big(\frac{1}{2}c\Big)}\nonumber\\
&\times(ch(\wedge^{2}(T_{\mathbf{C}}M)+\wedge^{2}(\xi_{\mathbf{C}})+T_{\mathbf{C}}M\otimes\xi_{\mathbf{C}}-15T_{\mathbf{C}}M-15\xi_{\mathbf{C}}+120)-\cosh\Big(\frac{1}{2}c\Big)ch(\wedge^{2}(T_{\mathbf{C}}M)\nonumber\\
&-2T_{\mathbf{C}}M\otimes\xi_{\mathbf{C}}+2\xi_{\mathbf{C}}\otimes\xi_{\mathbf{C}}-9T_{\mathbf{C}}M+22\xi_{\mathbf{C}}+39))\bigg\}^{(14)}=\bigg\{\widehat{A}(TM,\nabla^{TM})\frac{1}{2\sinh\Big(\frac{1}{2}c\Big)}\nonumber\\
&\times(ch(-\wedge^{3}(T_{\mathbf{C}}M)-\wedge^{3}(\xi_{\mathbf{C}})-T_{\mathbf{C}}M\otimes\wedge^{2}(\xi_{\mathbf{C}})-\xi_{\mathbf{C}}\otimes\wedge^{2}(T_{\mathbf{C}}M)+16\wedge^{2}(T_{\mathbf{C}}M)\nonumber\\
&+16\wedge^{2}(\xi_{\mathbf{C}})-T_{\mathbf{C}}M\otimes T_{\mathbf{C}}M+14T_{\mathbf{C}}M\otimes\xi_{\mathbf{C}}-\xi_{\mathbf{C}}\otimes\xi_{\mathbf{C}}-105T_{\mathbf{C}}M-105\xi_{\mathbf{C}}+567)\nonumber\\
&-\cosh\Big(\frac{1}{2}c\Big) ch(-\wedge^{3}(T_{\mathbf{C}}M)+\wedge^{3}(\xi_{\mathbf{C}})+S^{3}(\xi_{\mathbf{C}})-T_{\mathbf{C}}M\otimes\xi_{\mathbf{C}}\otimes\xi_{\mathbf{C}}-T_{\mathbf{C}}M\otimes\wedge^{2}(\xi_{\mathbf{C}})\nonumber\\
&-T_{\mathbf{C}}M\otimes S^{2}(\xi_{\mathbf{C}})+2\xi_{\mathbf{C}}\otimes\wedge^{2}(T_{\mathbf{C}}M)+\xi_{\mathbf{C}}\otimes\wedge^{2}(\xi_{\mathbf{C}})+\xi_{\mathbf{C}}\otimes S^{2}(\xi_{\mathbf{C}})+10\wedge^{2}(T_{\mathbf{C}}M)\nonumber\\
&-T_{\mathbf{C}}M\otimes T_{\mathbf{C}}M-20T_{\mathbf{C}}M\otimes\xi_{\mathbf{C}}+24\xi_{\mathbf{C}}\otimes\xi_{\mathbf{C}}-30T_{\mathbf{C}}M+100\xi_{\mathbf{C}}+70))\bigg\}^{(14)}.\nonumber
\end{align}
\end{thm}

\vskip 0.3 true cm
3.3. \noindent{ A general type of cancellation formulas for $2d$-dimensional Riemannian manifolds }
\vskip 0.3 true cm

Let $M$ be a $2d$-dimensional Riemannian manifold and set
\begin{align}
\Theta_{2}(T_{\mathbf{C}}M, m_{0}, \xi_{\mathbf{C}})=&\bigotimes^{\infty}_{n=1}S_{q^{n}}(\widetilde{T_{\mathbf{C}}M}-m_{0}\widetilde{\xi_{\mathbf{C}}})\otimes
\bigotimes_{m=1}^{\infty}\wedge_{-q^{m-\frac{1}{2}}}(\widetilde{T_{\mathbf{C}}M}-m_{0}\widetilde{\xi_{\mathbf{C}}}).
\end{align}
Clearly, $\Theta_{2}(T_{\mathbf{C}}M, m_{0}, \xi_{\mathbf{C}})$ admit formal Fourier expansion in $q^{\frac{1}{2}}$ as
\begin{align}
\Theta_{2}(T_{\mathbf{C}}M, m_{0}, \xi_{\mathbf{C}})=\widetilde{A_{0}}(T_{\mathbf{C}}M, m_{0}, \xi_{\mathbf{C}})+\widetilde{A_{1}}(T_{\mathbf{C}}M, m_{0}, \xi_{\mathbf{C}})q^{\frac{1}{2}}+\cdot\cdot\cdot,
\end{align}
where $\widetilde{A_{j}}$ are elements in the semi-group formally generated by Hermitian vector bundles over $M.$

Let $n$ be a nonnegative integer and satisfy $d-\left(2n+\frac{1-(-1)^d}{2}\right)>0$.
Choose
\begin{align}
\widetilde{P_{2}}(\tau)=&\bigg\{\widehat{A}(TM,\nabla^{TM})\Big(\sinh\Big(\frac{1}{2}c\Big)\Big)^{2n+\frac{1-(-1)^d}{2}}ch(\Theta_2(T_{\mathbf{C}}M, 2n+\frac{1-(-1)^d}{2}, \xi_{\mathbf{C}}))\bigg\}^{(2d)}.
\end{align}
We know $\widetilde{P_{2}}(\tau)$ is a modular form of weight $2m_{1}=d-\left(2n+\frac{1-(-1)^d}{2}\right)$ over $\Gamma^{0}(2),$ this can be found in \cite{HH,W1}.

Let us first observe that $m_{0}=2n+\frac{1-(-1)^d}{2}, m_{1}=1,$ then $d=2+m_{0}.$ Therefore
\begin{align}
\widetilde{P_{2}}(\tau)=&h_{0}(8\delta_{2}(\tau)).
\end{align}
Via simple calculations, we can obtain
\begin{thm}
For $(4+2m_{0})$-dimensional Riemannian manifold, we get
\begin{align}
&\bigg\{\widehat{A}(TM,\nabla^{TM})\Big(\sinh\Big(\frac{1}{2}c\Big)\Big)^{m_{0}}\bigg\}^{(4+2m_{0})}=\frac{1}{24}\bigg\{\widehat{A}(TM,\nabla^{TM})\Big(\sinh\Big(\frac{1}{2}c\Big)\Big)^{m_{0}}\\
&\times ch(-T_{\mathbf{C}}M+m_{0}\xi_{\mathbf{C}}+4)\bigg\}^{(4+2m_{0})}.\nonumber
\end{align}
\end{thm}
\begin{thm}
For $(4+2m_{0})$-dimensional Riemannian manifold, we can assert that
\begin{align}
&\bigg\{\widehat{A}(TM,\nabla^{TM})\Big(\sinh\Big(\frac{1}{2}c\Big)\Big)^{m_{0}}\bigg\}^{(4+2m_{0})}=\frac{1}{24}\bigg\{\widehat{A}(TM,\nabla^{TM})\Big(\sinh\Big(\frac{1}{2}c\Big)\Big)^{m_{0}}ch(\wedge^{2}(T_{\mathbf{C}}M)\\
&+m_{0}S^{2}(\xi_{\mathbf{C}})-m_{0}T_{\mathbf{C}}M\otimes\xi_{\mathbf{C}}+\frac{m_{0}(m_{0}-1)}{2}\xi_{\mathbf{C}}\otimes\xi_{\mathbf{C}}-3T_{\mathbf{C}}M+3m_{0}\xi_{\mathbf{C}}+6)\bigg\}^{(4+2m_{0})}.\nonumber
\end{align}
\end{thm}

The second part consider $m_{1}=2, d=4+m_{0}.$ Thus
\begin{align}
\widetilde{P_{2}}(\tau)=&h_{0}(8\delta_{2}(\tau))^{2}+h_{1}\varepsilon_{2}(\tau).
\end{align}
From this,
\begin{thm}
For $(8+2m_{0})$-dimensional Riemannian manifold, we see that
\begin{align}
&240\bigg\{\widehat{A}(TM,\nabla^{TM})\Big(\sinh\Big(\frac{1}{2}c\Big)\Big)^{m_{0}}\bigg\}^{(8+2m_{0})}+8\bigg\{\widehat{A}(TM,\nabla^{TM})\Big(\sinh\Big(\frac{1}{2}c\Big)\Big)^{m_{0}}ch(-T_{\mathbf{C}}M\\
&+m_{0}\xi_{\mathbf{C}}+8)\bigg\}^{(8+2m_{0})}=\bigg\{\widehat{A}(TM,\nabla^{TM})\Big(\sinh\Big(\frac{1}{2}c\Big)\Big)^{m_{0}}ch(\wedge^{2}(T_{\mathbf{C}}M)+m_{0}S^{2}(\xi_{\mathbf{C}})\nonumber\\
&-m_{0}T_{\mathbf{C}}M\otimes\xi_{\mathbf{C}}+\frac{m_{0}(m_{0}-1)}{2}\xi_{\mathbf{C}}\otimes\xi_{\mathbf{C}}-7T_{\mathbf{C}}M+7m_{0}\xi_{\mathbf{C}}+28)\bigg\}^{(8+2m_{0})}.\nonumber
\end{align}
\end{thm}

Assume that $m_{1}=3, d=6+m_{0}.$ it follows immediately that
\begin{thm}
For $(12+2m_{0})$-dimensional Riemannian manifold, we obtain
\begin{align}
&-504\bigg\{\widehat{A}(TM,\nabla^{TM})\Big(\sinh\Big(\frac{1}{2}c\Big)\Big)^{m_{0}}\bigg\}^{(12+2m_{0})}+32\bigg\{\widehat{A}(TM,\nabla^{TM})\Big(\sinh\Big(\frac{1}{2}c\Big)\Big)^{m_{0}}\\
&\times ch(-T_{\mathbf{C}}M+m_{0}\xi_{\mathbf{C}}+12)\bigg\}^{(12+2m_{0})}=\bigg\{\widehat{A}(TM,\nabla^{TM})\Big(\sinh\Big(\frac{1}{2}c\Big)\Big)^{m_{0}}ch(\wedge^{2}(T_{\mathbf{C}}M)\nonumber\\
&+m_{0}S^{2}(\xi_{\mathbf{C}})-m_{0}T_{\mathbf{C}}M\otimes\xi_{\mathbf{C}}+\frac{m_{0}(m_{0}-1)}{2}\xi_{\mathbf{C}}\otimes\xi_{\mathbf{C}}-11T_{\mathbf{C}}M+11m_{0}\xi_{\mathbf{C}}+66)\bigg\}^{(12+2m_{0})}.\nonumber
\end{align}
\end{thm}

In the following, let's expand $\Theta_{2}(T_{\mathbf{C}}M, m_{0}, \xi_{\mathbf{C}})$ when $m_{1}=4, d=8+m_{0}.$ In fact, we have
\begin{thm}
For $(16+2m_{0})$-dimensional Riemannian manifold, we deduce that
\begin{align}
&-7680\bigg\{\widehat{A}(TM,\nabla^{TM})\Big(\sinh\Big(\frac{1}{2}c\Big)\Big)^{m_{0}}\bigg\}^{(16+2m_{0})}+140\bigg\{\widehat{A}(TM,\nabla^{TM})\Big(\sinh\Big(\frac{1}{2}c\Big)\Big)^{m_{0}}\\
&\times ch(-T_{\mathbf{C}}M+m_{0}\xi_{\mathbf{C}}+16)\bigg\}^{(16+2m_{0})}+16\bigg\{\widehat{A}(TM,\nabla^{TM})\Big(\sinh\Big(\frac{1}{2}c\Big)\Big)^{m_{0}}ch(\wedge^{2}(T_{\mathbf{C}}M)\nonumber\\
&+m_{0}S^{2}(\xi_{\mathbf{C}})-m_{0}T_{\mathbf{C}}M\otimes\xi_{\mathbf{C}}+\frac{m_{0}(m_{0}-1)}{2}\xi_{\mathbf{C}}\otimes\xi_{\mathbf{C}}-15T_{\mathbf{C}}M+15m_{0}\xi_{\mathbf{C}}+120)\bigg\}^{(16+2m_{0})}\nonumber\\
&=\bigg\{\widehat{A}(TM,\nabla^{TM})\Big(\sinh\Big(\frac{1}{2}c\Big)\Big)^{m_{0}}ch(-\wedge^{3}(T_{\mathbf{C}}M)+m_{0}S^{3}(\xi_{\mathbf{C}})-\frac{m_{0}(m_{0}-1)}{2}T_{\mathbf{C}}M\otimes\xi_{\mathbf{C}}\nonumber\\
&\otimes\xi_{\mathbf{C}}-m_{0}T_{\mathbf{C}}M\otimes S^{2}(\xi_{\mathbf{C}})+m_{0}\xi_{\mathbf{C}}\otimes\wedge^{2}(T_{\mathbf{C}}M)+(m_{0}^{2}-m_{0})\xi_{\mathbf{C}}\otimes S^{2}(\xi_{\mathbf{C}})+16\wedge^{2}(T_{\mathbf{C}}M)\nonumber\\
&+16m_{0}S^{2}(\xi_{\mathbf{C}})-T_{\mathbf{C}}M\otimes T_{\mathbf{C}}M-14m_{0}T_{\mathbf{C}}M\otimes \xi_{\mathbf{C}}+(m_{0}^{3}+4m_{0}^{2}-6m_{0})\xi_{\mathbf{C}}\otimes\xi_{\mathbf{C}}\nonumber\\
&-105T_{\mathbf{C}}M+(-2m_{0}^{3}+6m_{0}^{2}+101m_{0})\xi_{\mathbf{C}}+\frac{4}{3}(m_{0}^{3}-3m_{0}^{2}+2m_{0}+432))\bigg\}^{(16+2m_{0})}.\nonumber
\end{align}
\end{thm}

\section{ $(a, b)$ type cancellation formulas }

\vskip 0.3 true cm
4.1. \noindent{ $(a, b)$ type cancellation formulas for $4d$-dimensional Riemannian manifolds }
\vskip 0.3 true cm

Set $M$ be a $4d$-dimensional Riemannian manifold with the associated Levi-Civita connection $\nabla^{TM}.$
Let $\xi$ be a complex line bundle and $\xi_{\mathbf{R}}$ be a rank two real oriented Euclidean vector bundle over $M$ carrying with a Euclidean connection $\nabla^{\xi_{\mathbf{R}}}$ which is the real bundle associated to $\xi.$
Write $a=(a_{1},\cdot\cdot\cdot, a_{k}), b=(b_{1},\cdot\cdot\cdot, b_{k}),$ where $a_{t}, b_{t}, 1\leqslant t\leqslant k$ are integers.

Suppose that
\begin{align}
\Theta_{2}(T_{\mathbf{C}}M, \xi_{\mathbf{R}}, a, b)=&\bigotimes^{\infty}_{n=1}S_{q^{n}}(\widetilde{T_{\mathbf{C}}M})\otimes\bigotimes^{\infty}_{m_{1}=1}\wedge_{q^{m_{1}}}(\widetilde{(\xi^{\otimes b_{1}})_{\mathbf{R}}\otimes\mathbf{C}})\otimes\cdot\cdot\cdot\\
&\otimes\bigotimes^{\infty}_{m_{k}=1}\wedge_{q^{m_{k}}}(\widetilde{(\xi^{\otimes b_{k}})_{\mathbf{R}}\otimes\mathbf{C}})\otimes\bigotimes^{\infty}_{r_{1}=1}\wedge_{q^{r_{1}-\frac{1}{2}}}(\widetilde{(\xi^{\otimes b_{1}})_{\mathbf{R}}\otimes\mathbf{C}})\otimes\cdot\cdot\cdot\nonumber\\
&\otimes\bigotimes^{\infty}_{r_{k}=1}\wedge_{q^{r_{k}-\frac{1}{2}}}(\widetilde{(\xi^{\otimes b_{k}})_{\mathbf{R}}\otimes\mathbf{C}})\otimes\bigotimes^{\infty}_{s_{1}=1}\wedge_{-q^{s_{1}-\frac{1}{2}}}(\widetilde{(\xi^{\otimes a_{1}})_{\mathbf{R}}\otimes\mathbf{C}})\otimes\cdot\cdot\cdot\nonumber\\
&\otimes\bigotimes^{\infty}_{s_{k}=1}\wedge_{-q^{s_{k}-\frac{1}{2}}}(\widetilde{(\xi^{\otimes a_{k}})_{\mathbf{R}}\otimes\mathbf{C}}).\nonumber
\end{align}
Obviously, $\Theta_{2}(T_{\mathbf{C}}M, \xi_{\mathbf{R}}, a, b)$ admit formal Fourier expansion in $q^{\frac{1}{2}}$ as
\begin{align}
\Theta_{2}(T_{\mathbf{C}}M, \xi_{\mathbf{R}}, a, b)=B_{0}(T_{\mathbf{C}}M, \xi_{\mathbf{R}}, a, b)+B_{1}(T_{\mathbf{C}}M, \xi_{\mathbf{R}}, a, b)q^{\frac{1}{2}}+\cdot\cdot\cdot,
\end{align}
where $B_{j}$ are elements in the semi-group formally generated by Hermitian vector bundles over $M.$

Let $\{ \pm2\pi\sqrt{-1}x_{j}|1\leq j\leq 2d\},$ $\{ \pm a_{t}2\pi\sqrt{-1}u\}$ and $\{ \pm b_{t}2\pi\sqrt{-1}u\}$ be the Chern roots of $T_{\mathbf{C}}M,$ $(\xi^{\otimes a_{t}})_{\mathbf{R}}\otimes\mathbf{C}$ and $(\xi^{\otimes b_{t}})_{\mathbf{R}}\otimes\mathbf{C}$ respectively.
Set $p_{1}$ denote the first Pontryagin form.
If $\omega$ is a differential form over $M,$ we denote by $\omega^{(i)}$ its degree $i$ component.
Write
\begin{align}
Q_{2}(\tau)=&\Big\{e^{\frac{1}{24}E_{2}(\tau)[p_{1}(TM)-\sum^{k}_{t=1}(a_{t}^{2}+2b_{t}^{2})p_{1}(\xi_{\mathbf{R}})]}2^{k}\cosh\Big(\frac{b_{1}}{2}c\Big)\cdot\cdot\cdot\cosh\Big(\frac{b_{k}}{2}c\Big)\\
&\times\widehat{A}(TM, \nabla^{TM})ch(\Theta_{2}(T_{\mathbf{C}}M, \xi_{\mathbf{R}}, a, b))\Big\}^{(4d)}.\nonumber
\end{align}
By \cite{LW2}, we find $Q_{2}(\tau)$ is a modular form of weight $2d$ over $\Gamma^{0}(2).$

Consider $d=2.$ Computations show that
\begin{thm}
For $8$-dimensional Riemannian manifold, we see that
\begin{align}
&[p_{1}(TM)-\sum^{k}_{t=1}(a_{t}^{2}+2b_{t}^{2})p_{1}(\xi_{\mathbf{R}})]\Big\{\frac{e^{\frac{1}{24}[p_{1}(TM)-\sum^{k}_{t=1}(a_{t}^{2}+2b_{t}^{2})p_{1}(\xi_{\mathbf{R}})]}+1}{p_{1}(TM)-\sum^{k}_{t=1}(a_{t}^{2}+2b_{t}^{2})p_{1}(\xi_{\mathbf{R}})}
2^{k}\cosh\Big(\frac{b_{1}}{2}c\Big)\cdot\cdot\cdot\\
&\times\cosh\Big(\frac{b_{k}}{2}c\Big)\widehat{A}(TM, \nabla^{TM})ch(240B_{0}^{1}(T_{\mathbf{C}}M, \xi_{\mathbf{R}}, a, b)+8B_{1}^{1}(T_{\mathbf{C}}M, \xi_{\mathbf{R}}, a, b)\nonumber\\
&-B_{2}^{1}(T_{\mathbf{C}}M, \xi_{\mathbf{R}}, a, b))+e^{\frac{1}{24}[p_{1}(TM)-\sum^{k}_{t=1}(a_{t}^{2}+2b_{t}^{2})p_{1}(\xi_{\mathbf{R}})]}2^{k}\cosh\Big(\frac{b_{1}}{2}c\Big)\cdot\cdot\cdot\cosh\Big(\frac{b_{k}}{2}c\Big)\nonumber\\
&\times\widehat{A}(TM, \nabla^{TM})ch(B_{0}^{1}(T_{\mathbf{C}}M, \xi_{\mathbf{R}}, a, b))\Big\}^{(4)}=\Big\{
2^{k}\cosh\Big(\frac{b_{1}}{2}c\Big) \cdot\cdot\cdot\cosh\Big(\frac{b_{k}}{2}c\Big)\nonumber\\
&\times \widehat{A}(TM, \nabla^{TM})ch(240B_{0}^{1}(T_{\mathbf{C}}M, \xi_{\mathbf{R}}, a, b)+8B_{1}^{1}(T_{\mathbf{C}}M, \xi_{\mathbf{R}}, a, b)-B_{2}^{1}(T_{\mathbf{C}}M, \xi_{\mathbf{R}}, a, b))\Big\}^{(8)},\nonumber
\end{align}
where
\begin{align}
B_{0}^{1}(T_{\mathbf{C}}M, \xi_{\mathbf{R}}, a, b)=&1;\\
B_{1}^{1}(T_{\mathbf{C}}M, \xi_{\mathbf{R}}, a, b)=&-((\xi^{\otimes a_{1}})_{\mathbf{R}}\otimes\mathbf{C}+\cdot\cdot\cdot+(\xi^{\otimes a_{k}})_{\mathbf{R}}\otimes\mathbf{C})\\
&+((\xi^{\otimes b_{1}})_{\mathbf{R}}\otimes\mathbf{C}+\cdot\cdot\cdot+(\xi^{\otimes b_{k}})_{\mathbf{R}}\otimes\mathbf{C});\nonumber\\
B_{2}^{1}(T_{\mathbf{C}}M, \xi_{\mathbf{R}}, a, b)=
&T_{\mathbf{C}}M-8-4k^{2}+4k+((\xi^{\otimes b_{1}})_{\mathbf{R}}\otimes\mathbf{C}+\cdot\cdot\cdot+(\xi^{\otimes b_{k}})_{\mathbf{R}}\otimes\mathbf{C})\\
&+(2k-2)((\xi^{\otimes a_{1}})_{\mathbf{R}}\otimes\mathbf{C}+\cdot\cdot\cdot+(\xi^{\otimes a_{k}})_{\mathbf{R}}\otimes\mathbf{C})\nonumber\\
&+(2k-2)((\xi^{\otimes b_{1}})_{\mathbf{R}}\otimes\mathbf{C}+\cdot\cdot\cdot+(\xi^{\otimes b_{k}})_{\mathbf{R}}\otimes\mathbf{C})\nonumber\\
&+[\wedge^{2}((\xi^{\otimes a_{1}})_{\mathbf{R}}\otimes\mathbf{C})+\cdot\cdot\cdot+\wedge^{2}((\xi^{\otimes a_{k}})_{\mathbf{R}}\otimes\mathbf{C})]\nonumber\\
&+[\wedge^{2}((\xi^{\otimes b_{1}})_{\mathbf{R}}\otimes\mathbf{C})+\cdot\cdot\cdot+\wedge^{2}((\xi^{\otimes b_{k}})_{\mathbf{R}}\otimes\mathbf{C})]\nonumber\\
&+[2(k-1)-((\xi^{\otimes a_{1}})_{\mathbf{R}}\otimes\mathbf{C}+\cdot\cdot\cdot+(\xi^{\otimes a_{k-1}})_{\mathbf{R}}\otimes\mathbf{C})]\nonumber\\
&\otimes(2-(\xi^{\otimes a_{k}})_{\mathbf{R}}\otimes\mathbf{C})+(2-(\xi^{\otimes a_{k-1}})_{\mathbf{R}}\otimes\mathbf{C})\nonumber\\
&\otimes[2(k-2)-((\xi^{\otimes a_{1}})_{\mathbf{R}}\otimes\mathbf{C}+\cdot\cdot\cdot+(\xi^{\otimes a_{k-2}})_{\mathbf{R}}\otimes\mathbf{C})]\nonumber\\
&+\cdot\cdot\cdot+(2-(\xi^{\otimes a_{1}})_{\mathbf{R}}\otimes\mathbf{C})\otimes(2-(\xi^{\otimes a_{2}})_{\mathbf{R}}\otimes\mathbf{C})\nonumber\\
&+[((\xi^{\otimes b_{1}})_{\mathbf{R}}\otimes\mathbf{C}+\cdot\cdot\cdot+(\xi^{\otimes b_{k-1}})_{\mathbf{R}}\otimes\mathbf{C})-2(k-1)]\nonumber\\
&\otimes((\xi^{\otimes b_{k}})_{\mathbf{R}}\otimes\mathbf{C}-2)+((\xi^{\otimes b_{k-1}})_{\mathbf{R}}\otimes\mathbf{C}-2)\nonumber\\
&\otimes[((\xi^{\otimes b_{1}})_{\mathbf{R}}\otimes\mathbf{C}+\cdot\cdot\cdot+(\xi^{\otimes b_{k-2}})_{\mathbf{R}}\otimes\mathbf{C})-2(k-2)]\nonumber\\
&+\cdot\cdot\cdot+((\xi^{\otimes b_{1}})_{\mathbf{R}}\otimes\mathbf{C}-2)\otimes((\xi^{\otimes b_{2}})_{\mathbf{R}}\otimes\mathbf{C}-2).\nonumber
\end{align}
\end{thm}

When $d=3,$ we can claim the following theorem.
\begin{thm}
For $12$-dimensional Riemannian manifold, we have
\begin{align}
&[p_{1}(TM)-\sum^{k}_{t=1}(a_{t}^{2}+2b_{t}^{2})p_{1}(\xi_{\mathbf{R}})]\Big\{\frac{e^{\frac{1}{24}[p_{1}(TM)-\sum^{k}_{t=1}(a_{t}^{2}+2b_{t}^{2})p_{1}(\xi_{\mathbf{R}})]}+1}{p_{1}(TM)-\sum^{k}_{t=1}(a_{t}^{2}+2b_{t}^{2})p_{1}(\xi_{\mathbf{R}})}
2^{k}\cosh\Big(\frac{b_{1}}{2}c\Big)\cdot\cdot\cdot\\
&\times\cosh\Big(\frac{b_{k}}{2}c\Big)\widehat{A}(TM, \nabla^{TM})ch(504B_{0}^{2}(T_{\mathbf{C}}M, \xi_{\mathbf{R}}, a, b)-32B_{1}^{2}(T_{\mathbf{C}}M, \xi_{\mathbf{R}}, a, b)\nonumber\\
&+B_{2}^{2}(T_{\mathbf{C}}M, \xi_{\mathbf{R}}, a, b))-e^{\frac{1}{24}[p_{1}(TM)-\sum^{k}_{t=1}(a_{t}^{2}+2b_{t}^{2})p_{1}(\xi_{\mathbf{R}})]}2^{k}\cosh\Big(\frac{b_{1}}{2}c\Big)\cdot\cdot\cdot\cosh\Big(\frac{b_{k}}{2}c\Big)\nonumber\\
&\times\widehat{A}(TM, \nabla^{TM})ch(B_{0}^{2}(T_{\mathbf{C}}M, \xi_{\mathbf{R}}, a, b))\Big\}^{(8)}=\Big\{
2^{k}\cosh\Big(\frac{b_{1}}{2}c\Big) \cdot\cdot\cdot\cosh\Big(\frac{b_{k}}{2}c\Big)\nonumber\\
&\times \widehat{A}(TM, \nabla^{TM})ch(504B_{0}^{2}(T_{\mathbf{C}}M, \xi_{\mathbf{R}}, a, b)-32B_{1}^{2}(T_{\mathbf{C}}M, \xi_{\mathbf{R}}, a, b)+B_{2}^{2}(T_{\mathbf{C}}M, \xi_{\mathbf{R}}, a, b))\Big\}^{(12)},\nonumber
\end{align}
where
\begin{align}
B_{0}^{2}(T_{\mathbf{C}}M, \xi_{\mathbf{R}}, a, b)=&1;\\
B_{1}^{2}(T_{\mathbf{C}}M, \xi_{\mathbf{R}}, a, b)=&-((\xi^{\otimes a_{1}})_{\mathbf{R}}\otimes\mathbf{C}+\cdot\cdot\cdot+(\xi^{\otimes a_{k}})_{\mathbf{R}}\otimes\mathbf{C})\\
&+((\xi^{\otimes b_{1}})_{\mathbf{R}}\otimes\mathbf{C}+\cdot\cdot\cdot+(\xi^{\otimes b_{k}})_{\mathbf{R}}\otimes\mathbf{C});\nonumber\\
B_{2}^{2}(T_{\mathbf{C}}M, \xi_{\mathbf{R}}, a, b)=
&T_{\mathbf{C}}M-12-4k^{2}+4k+((\xi^{\otimes b_{1}})_{\mathbf{R}}\otimes\mathbf{C}+\cdot\cdot\cdot+(\xi^{\otimes b_{k}})_{\mathbf{R}}\otimes\mathbf{C})\\
&+(2k-2)((\xi^{\otimes a_{1}})_{\mathbf{R}}\otimes\mathbf{C}+\cdot\cdot\cdot+(\xi^{\otimes a_{k}})_{\mathbf{R}}\otimes\mathbf{C})\nonumber\\
&+(2k-2)((\xi^{\otimes b_{1}})_{\mathbf{R}}\otimes\mathbf{C}+\cdot\cdot\cdot+(\xi^{\otimes b_{k}})_{\mathbf{R}}\otimes\mathbf{C})\nonumber\\
&+[\wedge^{2}((\xi^{\otimes a_{1}})_{\mathbf{R}}\otimes\mathbf{C})+\cdot\cdot\cdot+\wedge^{2}((\xi^{\otimes a_{k}}_{\mathbf{R}})\otimes\mathbf{C})]\nonumber\\
&+[\wedge^{2}((\xi^{\otimes b_{1}})_{\mathbf{R}}\otimes\mathbf{C})+\cdot\cdot\cdot+\wedge^{2}((\xi^{\otimes b_{k}}_{\mathbf{R}})\otimes\mathbf{C})]\nonumber\\
&+[2(k-1)-((\xi^{\otimes a_{1}})_{\mathbf{R}}\otimes\mathbf{C}+\cdot\cdot\cdot+(\xi^{\otimes a_{k-1}})_{\mathbf{R}}\otimes\mathbf{C})]\nonumber\\
&\otimes(2-(\xi^{\otimes a_{k}})_{\mathbf{R}}\otimes\mathbf{C})+(2-(\xi^{\otimes a_{k-1}})_{\mathbf{R}}\otimes\mathbf{C})\nonumber\\
&\otimes[2(k-2)-((\xi^{\otimes a_{1}})_{\mathbf{R}}\otimes\mathbf{C}+\cdot\cdot\cdot+(\xi^{\otimes a_{k-2}})_{\mathbf{R}}\otimes\mathbf{C})]\nonumber\\
&+\cdot\cdot\cdot+(2-(\xi^{\otimes a_{1}})_{\mathbf{R}}\otimes\mathbf{C})\otimes(2-(\xi^{\otimes a_{2}})_{\mathbf{R}}\otimes\mathbf{C})\nonumber\\
&+[((\xi^{\otimes b_{1}})_{\mathbf{R}}\otimes\mathbf{C}+\cdot\cdot\cdot+(\xi^{\otimes b_{k-1}})_{\mathbf{R}}\otimes\mathbf{C})-2(k-1)]\nonumber\\
&\otimes((\xi^{\otimes b_{k}})_{\mathbf{R}}\otimes\mathbf{C}-2)+((\xi^{\otimes b_{k-1}})_{\mathbf{R}}\otimes\mathbf{C}-2)\nonumber\\
&\otimes[((\xi^{\otimes b_{1}})_{\mathbf{R}}\otimes\mathbf{C}+\cdot\cdot\cdot+(\xi^{\otimes b_{k-2}})_{\mathbf{R}}\otimes\mathbf{C})-2(k-2)]\nonumber\\
&+\cdot\cdot\cdot+((\xi^{\otimes b_{1}})_{\mathbf{R}}\otimes\mathbf{C}-2)\otimes((\xi^{\otimes b_{2}})_{\mathbf{R}}\otimes\mathbf{C}-2).\nonumber
\end{align}
\end{thm}

\vskip 0.3 true cm
4.2. \noindent{ $(a, b)$ type cancellation formulas involving a complex line bundle for $4d$-dimensional Riemannian manifolds }
\vskip 0.3 true cm

Let $\eta$ be a rank two oriented Euclidean vector bundle over $M$ carrying a Euclidean connection $\nabla^{\eta}.$
Assume that
\begin{align}
\overline{\Theta_{2}}(T_{\mathbf{C}}M, \xi_{\mathbf{R}}, \eta_{\mathbf{C}}, a, b)=&\bigotimes^{\infty}_{n=1}S_{q^{n}}(\widetilde{T_{\mathbf{C}}M})\otimes\bigotimes^{\infty}_{m=1}\wedge_{q^{m}}(\widetilde{(\xi^{\otimes b})_{\mathbf{R}}\otimes\mathbf{C}}+\widetilde{\eta_{\mathbf{C}}})\\
&\otimes\bigotimes^{\infty}_{r=1}\wedge_{q^{r-\frac{1}{2}}}(\widetilde{(\xi^{\otimes b})_{\mathbf{R}}\otimes\mathbf{C}}+\widetilde{\eta_{\mathbf{C}}})\otimes\bigotimes^{\infty}_{s=1}\wedge_{-q^{s-\frac{1}{2}}}(\widetilde{(\xi^{\otimes a})_{\mathbf{R}}\otimes\mathbf{C}}-2\widetilde{\eta_{\mathbf{C}}}).\nonumber
\end{align}

In the meantime we have
\begin{align}
\overline{Q_{2}}(\tau)=&\Big\{e^{\frac{1}{24}E_{2}(\tau)[p_{1}(TM)-(a^{2}+2b^{2})p_{1}(\xi_{\mathbf{R}})]}2\cosh\Big(\frac{b}{2}c\Big)\cosh\Big(\frac{1}{2}\overline{c}\Big)\\
&\times\widehat{A}(TM, \nabla^{TM})ch(\overline{\Theta_{2}}(T_{\mathbf{C}}M, \xi_{\mathbf{R}}, \eta_{\mathbf{C}}, a, b))\Big\}^{(4d)},\nonumber
\end{align}
where $\overline{c}=2\pi\sqrt{-1}\overline{u}.$
In \cite{LW2}, it is shown that $\overline{Q_{2}}(\tau)$ is a modular form of weight $2d$ over $\Gamma^{0}(2).$

We consider $d=2.$ In this case, we have
\begin{thm}
For $8$-dimensional Riemannian manifold, we conclude that
\begin{align}
&[p_{1}(TM)-(a^{2}+2b^{2})p_{1}(\xi_{\mathbf{R}})]\Big\{\frac{e^{\frac{1}{24}[p_{1}(TM)-(a^{2}+2b^{2})p_{1}(\xi_{\mathbf{R}})]}+1}{p_{1}(TM)-(a^{2}+2b^{2})p_{1}(\xi_{\mathbf{R}})}
2\cosh\Big(\frac{b}{2}c\Big)\cosh\Big(\frac{1}{2}\overline{c}\Big)\\
&\times\widehat{A}(TM, \nabla^{TM})ch(240\overline{B_{0}^{1}}(T_{\mathbf{C}}M, \xi_{\mathbf{R}}, \eta_{\mathbf{C}}, a, b)+8\overline{B_{1}^{1}}(T_{\mathbf{C}}M, \xi_{\mathbf{R}}, \eta_{\mathbf{C}}, a, b)\nonumber\\
&-\overline{B_{2}^{1}}(T_{\mathbf{C}}M, \xi_{\mathbf{R}}, \eta_{\mathbf{C}}, a, b))+e^{\frac{1}{24}[p_{1}(TM)-(a^{2}+2b^{2})p_{1}(\xi_{\mathbf{R}})]}2\cosh\Big(\frac{b}{2}c\Big)\cosh\Big(\frac{1}{2}\overline{c}\Big)\nonumber\\
&\times\widehat{A}(TM, \nabla^{TM})ch(\overline{B_{0}^{1}}(T_{\mathbf{C}}M, \xi_{\mathbf{R}}, \eta_{\mathbf{C}}, a, b))\Big\}^{(4)}=\Big\{
2\cosh\Big(\frac{b}{2}c\Big)\cosh\Big(\frac{1}{2}\overline{c}\Big)\nonumber\\
&\times \widehat{A}(TM, \nabla^{TM})ch(240\overline{B_{0}^{1}}(T_{\mathbf{C}}M, \xi_{\mathbf{R}}, \eta_{\mathbf{C}}, a, b)+8\overline{B_{1}^{1}}(T_{\mathbf{C}}M, \xi_{\mathbf{R}}, \eta_{\mathbf{C}}, a, b)\nonumber\\
&-\overline{B_{2}^{1}}(T_{\mathbf{C}}M, \xi_{\mathbf{R}}, \eta_{\mathbf{C}}, a, b))\Big\}^{(8)},\nonumber
\end{align}
where
\begin{align}
\overline{B_{0}^{1}}(T_{\mathbf{C}}M, \xi_{\mathbf{R}}, \eta_{\mathbf{C}}, a, b)=&1;\\
\overline{B_{1}^{1}}(T_{\mathbf{C}}M, \xi_{\mathbf{R}}, \eta_{\mathbf{C}}, a, b)=&-(\xi^{\otimes a})_{\mathbf{R}}\otimes\mathbf{C}+(\xi^{\otimes b})_{\mathbf{R}}\otimes\mathbf{C}+3\eta_{\mathbf{C}}-6;\\
\overline{B_{2}^{1}}(T_{\mathbf{C}}M, \xi_{\mathbf{R}}, \eta_{\mathbf{C}}, a, b)=
&\wedge^{2}((\xi^{\otimes a})_{\mathbf{R}}\otimes\mathbf{C})+\wedge^{2}((\xi^{\otimes b})_{\mathbf{R}}\otimes\mathbf{C})-((\xi^{\otimes a})_{\mathbf{R}}\otimes\mathbf{C})\\
&\otimes((\xi^{\otimes b})_{\mathbf{R}}\otimes\mathbf{C})-3((\xi^{\otimes a})_{\mathbf{R}}\otimes\mathbf{C})\otimes\eta_{\mathbf{C}}\nonumber\\
&+3((\xi^{\otimes b})_{\mathbf{R}}\otimes\mathbf{C})\otimes\eta_{\mathbf{C}}+4\eta_{\mathbf{C}}\otimes\eta_{\mathbf{C}}+S^{2}(\eta_{\mathbf{C}})+T_{\mathbf{C}}M\nonumber\\
&+6(\xi^{\otimes a})_{\mathbf{R}}\otimes\mathbf{C}-5(\xi^{\otimes b})_{\mathbf{R}}\otimes\mathbf{C}-17\eta_{\mathbf{C}}+7.\nonumber
\end{align}
\end{thm}

Fix $d=3,$ it is easily seen that
\begin{thm}
For $12$-dimensional Riemannian manifold, we can assert that
\begin{align}
&[p_{1}(TM)-(a^{2}+2b^{2})p_{1}(\xi_{\mathbf{R}})]\Big\{\frac{e^{\frac{1}{24}[p_{1}(TM)-(a^{2}+2b^{2})p_{1}(\xi_{\mathbf{R}})]}+1}{p_{1}(TM)-(a^{2}+2b^{2})p_{1}(\xi_{\mathbf{R}})}
2\cosh\Big(\frac{b}{2}c\Big)\cosh\Big(\frac{1}{2}\overline{c}\Big)\\
&\times\widehat{A}(TM, \nabla^{TM})ch(504\overline{B_{0}^{2}}(T_{\mathbf{C}}M, \xi_{\mathbf{R}}, \eta_{\mathbf{C}}, a, b)-32\overline{B_{1}^{2}}(T_{\mathbf{C}}M, \xi_{\mathbf{R}}, \eta_{\mathbf{C}}, a, b)\nonumber\\
&+\overline{B_{2}^{2}}(T_{\mathbf{C}}M, \xi_{\mathbf{R}}, \eta_{\mathbf{C}}, a, b))-e^{\frac{1}{24}[p_{1}(TM)-(a^{2}+2b^{2})p_{1}(\xi_{\mathbf{R}})]}2\cosh\Big(\frac{b}{2}c\Big)\cosh\Big(\frac{1}{2}\overline{c}\Big)\nonumber\\
&\times\widehat{A}(TM, \nabla^{TM})ch(\overline{B_{0}^{2}}(T_{\mathbf{C}}M, \xi_{\mathbf{R}}, \eta_{\mathbf{C}}, a, b))\Big\}^{(8)}=\Big\{
2\cosh\Big(\frac{b}{2}c\Big)\cosh\Big(\frac{1}{2}\overline{c}\Big)\nonumber\\
&\times \widehat{A}(TM, \nabla^{TM})ch(504\overline{B_{0}^{2}}(T_{\mathbf{C}}M, \xi_{\mathbf{R}}, \eta_{\mathbf{C}}, a, b)-32\overline{B_{1}^{2}}(T_{\mathbf{C}}M, \xi_{\mathbf{R}}, \eta_{\mathbf{C}}, a, b)\nonumber\\
&+\overline{B_{2}^{2}}(T_{\mathbf{C}}M, \xi_{\mathbf{R}}, \eta_{\mathbf{C}}, a, b))\Big\}^{(12)},\nonumber
\end{align}
where
\begin{align}
\overline{B_{0}^{2}}(T_{\mathbf{C}}M, \xi_{\mathbf{R}}, \eta_{\mathbf{C}}, a, b)=&1;\\
\overline{B_{1}^{2}}(T_{\mathbf{C}}M, \xi_{\mathbf{R}}, \eta_{\mathbf{C}}, a, b)=&-(\xi^{\otimes a})_{\mathbf{R}}\otimes\mathbf{C}+(\xi^{\otimes b})_{\mathbf{R}}\otimes\mathbf{C}+3\eta_{\mathbf{C}}-6;\\
\overline{B_{2}^{2}}(T_{\mathbf{C}}M, \xi_{\mathbf{R}}, \eta_{\mathbf{C}}, a, b)=
&\wedge^{2}((\xi^{\otimes a})_{\mathbf{R}}\otimes\mathbf{C})+\wedge^{2}((\xi^{\otimes b})_{\mathbf{R}}\otimes\mathbf{C})-((\xi^{\otimes a})_{\mathbf{R}}\otimes\mathbf{C})\\
&\otimes((\xi^{\otimes b})_{\mathbf{R}}\otimes\mathbf{C})-3((\xi^{\otimes a})_{\mathbf{R}}\otimes\mathbf{C})\otimes\eta_{\mathbf{C}}\nonumber\\
&+3((\xi^{\otimes b})_{\mathbf{R}}\otimes\mathbf{C})\otimes\eta_{\mathbf{C}}+4\eta_{\mathbf{C}}\otimes\eta_{\mathbf{C}}+S^{2}(\eta_{\mathbf{C}})+T_{\mathbf{C}}M\nonumber\\
&+6(\xi^{\otimes a})_{\mathbf{R}}\otimes\mathbf{C}-5(\xi^{\otimes b})_{\mathbf{R}}\otimes\mathbf{C}-17\eta_{\mathbf{C}}+3.\nonumber
\end{align}
\end{thm}

\section{ Acknowledgements}

The author was supported in part by  NSFC No.11771070. The author thanks the referee for his (or her) careful reading and helpful comments.

\vskip 1 true cm


\bigskip
\bigskip

\noindent {\footnotesize {\it S. Liu} \\
{School of Mathematics and Statistics, Northeast Normal University, Changchun 130024, China}\\
{Email: liusy719@nenu.edu.cn}

\noindent {\footnotesize {\it Y. Wang} \\
{School of Mathematics and Statistics, Northeast Normal University, Changchun 130024, China}\\
{Email: wangy581@nenu.edu.cn}

\clearpage
\section*{Statements and Declarations}

Funding: This research was funded by National Natural Science Foundation of China: No.11771070.\\

Competing Interests: The authors have no relevant financial or non-financial interests to disclose.\\

Author Contributions: All authors contributed to the study conception and design. Material preparation, data collection and analysis were performed by Siyao Liu and Yong Wang. The first draft of the manuscript was written by Siyao Liu and all authors commented on previous versions of the manuscript. All authors read and approved the final manuscript.\\

Availability of Data and Material: The datasets supporting the conclusions of this article are included within the article and its additional files.\\

\end{document}